\documentclass[reqno, 11pt]{amsart} \usepackage{various, geometry}
\usepackage{hyperref}

\begin{document}

\title{Squarefree numbers in short intervals: explicit and formalized}
\author{Mayank Pandey}
\maketitle

\section{Introduction}

In this note, we give an account of explicating the exponent (and along the way, all
implicit constants) of the result in \cite{2024arXiv2401.13981} on squarefree numbers in short intervals, while at
the same time formalizing the proof in Lean 4
\footnote{Formalization located in this repository: \url{https://github.com/mayankpandey9973/squarefree-explicit-formalized}\raggedright}. The formalization was performed by a combination of GPT 5.5 and Claude Opus 4.8.

The explicit and formalized result is as follows.
\begin{theorem}\label{theorem:cuo7a3hjxd}
  Suppose that $0 < \varepsilon\le 1/90935 $ and that
  $X\ge \exp(10^{27}/\varepsilon^2) $. Then, with $H =  X^{1/5 - 2/90935 + \varepsilon}$, we have that
  \begin{equation}
    \biggl|\sum_{X\le n\le X + H } \mu(n)^2 - \frac{6}{\pi^2}H\biggr|
    \le \frac{10^{450}}{\varepsilon} H X^{-\varepsilon/10^{25}}.
  \end{equation}
\end{theorem}
This is, in the repository, located at
\nolinkurl{squarefree_lean/Squarefree/ExplicitBounds.lean:765}, as the theorem
\nolinkurl{count_short_interval_eps}.

The purpose of this note is to give a rough outline of the additions mathematically that
went into making the exponent explicit (\cite{2024arXiv2401.13981} depends on work of
Green and Tao \cite{GT} on the quantitative equidistribution of nilsequences which does not
have explicit exponents).

\section{Layout of argument}

We give in this section a rough outline of the argument so that we may make clear where
the work in making the exponent explicit went. We refer readers to \cite[\S2]{2024arXiv2401.13981} for a more
complete account of the setup. To simplify, we'll state bounds that yield an exponent of
$2/92485$. See Remark \ref{remark:cuo6920uql} for more details of the sacrifice made. 

To estimate (with a power savings, say) the quantity
\begin{equation}
  \sum_{X\le n\le X + H } \mu(n)^2, 
\end{equation}
we are reduced to showing bounds of $O(H X^{-\varepsilon}) $ on the size of 
\begin{equation}
  N(H, D)  := \#\set{d\in [D, 2D] : \| X/d^2 \|\le H/D^2}
\end{equation}
for all $H X^{-\varepsilon}\le D\le\sqrt{X} $. To explain the full proof of Theorem \ref{theorem:cuo7a3hjxd},
we must bound a somewhat finer quantity, but we weaken the result slightly here for clarity.
See Remark \ref{remark:cuo6920uql} for further detail.

If $H $ is a fixed power of $X $, obtaining such bounds for $D\ll H X^{\delta} $
for some $\delta > 0 $ is standard, either by repeated differencing or standard exponential sum
techniques. 

This is achieved with the following bound (for convenience, we will from now on be writing
$\Delta := D/H$, $G := X/H^5 $ and we'll suppose that $G\ge 1 $) and simply follows from repeated differencing:
\begin{proposition}[Via \nolinkurl{Squarefree.Counting.fourthDeriv_count}, as used in the proof of \nolinkurl{Squarefree.prop_2_4}]\label{proposition:cun54m4irw}
  \begin{equation}
    N(H, D) \ll H \biggl( \frac{\Delta^{7/8}}{H^{1/8}} + \frac{\Delta^{3/4}}{H^{1/8}} + \frac{\Delta^{11/8}}{G^{1/8}H^{1/8}}
    + \frac{G^{1/15}\Delta^{3/5}}{H^{1/15}}\biggr).
  \end{equation}
\end{proposition}
We have omitted explicit constants for the sake of readability.

This role was played in the original paper by \cite[Proposition 2.2]{2024arXiv2401.13981}, which
is a bit stronger then \ref{proposition:cun54m4irw}, but Proposition \ref{proposition:cun54m4irw} suffices for our purposes.

For the regimes $\Delta\gg H^{1/2} X^{\varepsilon} $ and $\Delta\ll H^{1/2} X^{-\varepsilon} $, we have Propositions \ref{proposition:cun54s0uod} and \ref{proposition:cun54s2ei7}, respectively, below.
These correspond to (4) and (3) of \cite{2024arXiv2401.13981}, respectively.
\begin{proposition}[Via \nolinkurl{Squarefree.prop_3_2_fiber}, \nolinkurl{Squarefree.prop_6_1}]\label{proposition:cun54s0uod}
With $C = 10^{20} $, for any $\log X \le  U\le X^{1/100}$, we have that
\begin{equation}
  N(H, D) \ll H\max \biggl( \frac{1}{U}, U^C \biggl( \frac{H^{2/3}}{\Delta^{4/3}}G^{4/3} + \frac{H}{\Delta^2}G
  + \frac{1}{H^{1/2}}G^{1/2}\biggr)  \biggr)
\end{equation}
\end{proposition}

\begin{proposition}[Via \nolinkurl{Squarefree.prop_5_1}]\label{proposition:cun54s2ei7}
  Suppose that $\log X \le  U\le X^{1/100}$, $\Delta \le H^{1/2} G^{-1/2}U^{-5}$, and $\Delta\gg G^2U^5 + G^4U^{20} $.
  Then, we have 
  \begin{equation}
    N(H, D)\ll H \max\biggl( \frac{1}{U}, \frac{1}{\Delta^{1/2}}G^9U^{51} + \frac{1}{\Delta}G^{19}U^{96} + \frac{\Delta^2}{H}G^{45/2}U^{243/2}  \biggr).
  \end{equation}
\end{proposition}
Proposition \ref{proposition:cun54s0uod} is slightly stronger than the corresponding \cite[(4)]{2024arXiv2401.13981} owing to our decision to use
elementary bounds due to Swinnerton-Dyer \cite{SWINNERTONDYER1974128} to avoid dealing with any Fourier analysis to simplify the
formalization. In the original paper, the $H^{2/3}/\Delta^{4/3} $ term would be weakened to, up to a factor of $(GU)^{O(1)} $,
$H^{1/2}/\Delta $. The bottleneck for the exponent does not involve Proposition \ref{proposition:cun54s0uod}, nor would
it if the weaker \cite[(4)]{2024arXiv2401.13981} were used instead.

Proposition \ref{proposition:cun54s2ei7} is essentially \cite[(3)]{2024arXiv2401.13981}, though with slightly better exponents, not
due to distinct arguments but to a slightly more careful gluing-together of different ranges. 
Together, Propositions \ref{proposition:cun54m4irw}, \ref{proposition:cun54s0uod}, and \ref{proposition:cun54s2ei7} imply that $N(H, D) \ll H X^{-\varepsilon}$
for all ranges besides
\begin{equation}\label{eq:cuo2oevxbi}
  H^{1/2} G^{-45/4} U^{-O(1)} \ll\Delta\ll H^{1/2} GU^{O(1)}  .
\end{equation}
This is treated by \cite[(5)]{2024arXiv2401.13981}, which is ultimately shown with an appeal to \cite{GT}. As a result, the
exponent obtained is not explicit. In this coming section, we sketch the observations that allow us to,
without making all of the rather general result in \cite{GT} explicit (which would likely have either been
much more cumbersome and led to a much weaker exponent), deal with our specific situation manually. It is worth noting that we are able to succeed here due to the underlying
polynomial orbit on a nilpotent Lie group present (at least after looking at a short interval), so we
are not straying far from the ideas of \cite{GT}.

\section{Explicit bracket expression nonconcentration}
We are roughly reduced, in the critical case that remains (see \cite[\S5.2]{2024arXiv2401.13981} for the
reduction in the original paper), to showing that 
\begin{equation}\label{eq:cuo2oc404s}
 g(r) = f_3(r) + f_1(r) \{f_2(r)\}
\end{equation}
cannot be $O(X^{-1/10}) \mod 1$ for a $\gg X^{-o(1)} $-fraction of $r\sim X^{1/10} $, where $f_1, f_2, f_3$ are,
roughly, monomials of size $X^{1/10}, X^{1/5}, X^{3/10}$. This is resolved in \cite{2024arXiv2401.13981} by splitting $r $ into short
intervals, Taylor expanding the $f_i $, applying the main theorem of \cite{GT}, and showing simply
that the corresponding polynomial orbit on a nilmanifold equidistributes (at scales of size $X^{-\delta} $)
for most short intervals.

We instead carry this out bare-handedly, applying van der Corput differencing directly to the bracket expression
\eqref{eq:cuo2oc404s}. Since we only care about upper bounds, rather than equidistribution, our job is simplified in two ways.
First, we can, directly, assuming the negation, suppose that there exist $h_1, h_2, h_3\ll X^{O(\delta)} $ (say, from the start, we
seek a small power saving) for which $\partial_{h_1, h_2, h_3}g(r) $
is $O(X^{-1/10})$ for a $\gg X^{-o(1)} $-fraction of integers. Here, we write $\partial_h g(r) = g(r) - g(r + h)  $ and $\partial_{h_1,\dots,h_j} = \partial_{h_1}\dots\partial_{h_j} $. 
The second major simplification comes from the fact that we do not need to detect wraparounds with the fractional
part: since
$\{x\} - \{y\}\in \{\{x - y\}, \{x - y\} - 1\} $, we can simply, at the cost of splitting over finitely many cases, fiber over these
possibilities. More precisely, we note that 
\begin{multline}\label{eq:cuo2od10w7}
  \partial_{h_1,h_2,h_3}g(r) = \partial_{h_1,h_2,h_3}f_3(r) + f_1(r+h_1 + h_2 + h_3)(\partial_{h_1,h_2,h_3}{f_2}(r) + \rho_0) \\  + \sum_{i\le 3 \\ j < k\text{ s.t. }\{1,2,3\} = \{i,j,k\} }^3
  \partial_{h_i}f_1(r+ h_1 + h_2 + h_3 - h_i)(\{\partial_{h_j,h_k}{f_2}(r)\} + \rho_i) + O(X^{-1/10 + O(\delta)})
\end{multline}
where $\rho_0, \rho_1, \rho_2, \rho_3 $ are bounded integers carrying the information of carrying (for example,
$\rho_0$ corresponds to the value of $ \partial_{h_1, h_2, h_3}\floor{f_2}\in \{-5, \dots, 5\}  $ and the other $\rho_i $ are similarly defined
integers of size $O(1) $). 

Rather than trying to detect the $\rho_i$ as would work to show equidistribution, we simply split
into cases with a union bound. The quantity within the fractional parts in \eqref{eq:cuo2od10w7} is always $O(X^{O(\delta)}) $,
so we may, at the cost of $X^{O(\delta)} $ possible integer shifts, remove the fractional part. This
sum of (effectively up to $2 $, it turns out) monomials can then be shown to be  $O(X^{-1/10 + O(\delta)}) $
rarely by standard elementary techniques (those of \cite[\S3]{2024arXiv2401.13981} suffice, it turns out, though we
use the slightly better Swinnerton-Dyer estimate from \cite{SWINNERTONDYER1974128} in the formalization in place of
\cite[Proposition 3.3]{2024arXiv2401.13981}). 

In the end, we show the following bound:
\begin{proposition}[Unoptimized full range consequence of \nolinkurl{Squarefree.prop_7_3_explicit}, as assembled in \nolinkurl{Squarefree.dblock_on_strip_explicit}]\label{proposition:cuo69xea7p}
  Write
  \begin{equation}
    \Psi_{e} = \biggl(\frac{\Delta^2}{H}\biggr)^e + \biggl(\frac{\Delta^2}{H}\biggr)^{e - 1}.
  \end{equation}
  We have that for $G\le H^{1/10}  $, 
  \begin{align*}
    N(H, D)\ll H \biggl( \frac{1}{U} &+ 
    U^{O(1)} (H^{-1/84}G^{6/7}\Psi_{5/84}
                                     + H^{-1/42}G^{22/21}\Psi_{1/42} + H^{-1/16}G\Psi_{5/16} \\
    &+ H^{-1/28}G^{13/14}\Psi_{5/28} 
      + H^{-1/18}G^{4/3}\Psi_{1/18} + H^{-1/54}G^{25/27}\Psi_{1/54}  \\
    &+ H^{-1/28}G^{13/14}\Psi_{-1/28}
    + H^{-1/6}G^{11/9}\Psi_{-1/6} + H^{-1/16}G\Psi_{-1/16}) \biggr).
  \end{align*}

\end{proposition}

\begin{remark}\label{remark:cuo6920uql}
  With both Propositions \ref{proposition:cun54s2ei7} and \ref{proposition:cuo69xea7p}, we have weakened (in a fashion that affects the
  final exponent) the results shown in the Lean formalization for
  simplicity of exposition. In reality, we bound the size of
  \begin{equation*}
    N_a(H, D) = \#\left\{d\in [D, 2D] : \left\| \frac{X}{d^2} \right\|, \left\| \frac{X}{(d + a)^2} \right\| \le \frac{H}{D^2} \land
      \left\| \frac{X}{d'^2} \right\| > \frac{H}{D^2}\forall d< d' < d + a \right\}.
  \end{equation*}
  Initial reductions show that if $N(H, D)\gg H/U $, it suffices to deal with
  \begin{equation}
    \Delta G^{-1/4}U^{-3/4}\ll a\ll \Delta U,
  \end{equation}
  which is a fairly narrow region. 
  Nevertheless, around the critical region, whether we apply \nolinkurl{Squarefree.prop_5_1} or
  \nolinkurl{Squarefree.prop_7_3} depends not only on the scale of $\Delta $ but on that of $a $.
  
  We have elected to give slightly less precise statements for we do not feel there is
  sufficient mathematical content within this slight optimization to warrant complicating
  the existing note.
\end{remark}

\bibliography{refs}{} \bibliographystyle{alpha}
\end{document}